\input amstex
 \documentstyle{amsppt}
 \magnification=\magstep1
 \vsize=24.2true cm
 \hsize=15.3true cm
 \nopagenumbers\topskip=1truecm
 \headline={\tenrm\hfil\folio\hfil}

 \TagsOnRight

\hyphenation{auto-mor-phism auto-mor-phisms co-homo-log-i-cal co-homo-logy
co-homo-logous dual-izing pre-dual-izing geo-metric geo-metries geo-metry
half-space homeo-mor-phic homeo-mor-phism homeo-mor-phisms homo-log-i-cal
homo-logy homo-logous homo-mor-phism homo-mor-phisms hyper-plane hyper-planes
hyper-sur-face hyper-sur-faces idem-potent iso-mor-phism iso-mor-phisms
multi-plic-a-tion nil-potent poly-nomial priori rami-fication sin-gu-lar-ities
sub-vari-eties sub-vari-ety trans-form-a-tion trans-form-a-tions Castel-nuovo
Enri-ques Lo-ba-chev-sky Theo-rem Za-ni-chelli in-vo-lu-tion Na-ra-sim-han Bohr-Som-mer-feld}

\define\rest#1{_{\textstyle{|}#1}} 

\define\Span#1{\left<#1\right>} 

\define\half{{\textstyle{1\over2}}}
\define\wave{\widetilde}

\define\C{\Bbb C} 
\define\R{\Bbb R} 
\define\Z{\Bbb Z} 

\define\proj{\Bbb P} 

\define\sA{{\Cal A}} 
\define\sG{{\Cal G}} 
\define\sM{{\Cal M}} 

\define\al{\alpha}
\define\be{\beta}

\define\om{\omega}
\define\si{\sigma}

\define\Ga{\Gamma}
\define\La{\Lambda}
\define\Om{\Omega}









 \document

  \topmatter
  \title Pseudo symplectic geometry as an extension of the
symplectic geometry \endtitle
  \author Nik. Tyurin                   \endauthor

   \address MPI, Bonn
   \endaddress
  \email
    tyurin\@tyurin.mccme.ru  jtyurin\@mpim-bonn.mpg.de
   \endemail

\abstract In [9] we've outlined how the Taubes result about nontriviality
of the Seiberg --- Witten invariants for symplectic manifolds can be gene ralized to more general case, when a priori based manifold doesn't admit a    symplectic structure.
In the present paper we correct our result from [9] and discuss some properties which belong to pseudo symplectic
4- manioflds in the framework of the complex geometry.
 
\endabstract

   \endtopmatter

\head \S 0. Introduction  \endhead                                               
$$
$$

In 1994 new invariants of smooth structures in 4- dimensional geometry
were introduced by E.Witten in his seminal work [10]. These invariants are very closed to the Donaldson invariants (and a relationship was established in [5]).
But so as it's much more easy to compute the Seiberg --- Witten invariants
a number of unsolved problems have been solved in the framework of the 
Seiberg --- Witten theory.  Two first results proved in [10] have been proved
earlier in the framework of the Donaldson theory, namely

1. that for a connected sum $Y = X_1 \sharp X_2$ where $b_2^+(X_i) >0$
 the invariants vanish (as well as for the Donaldson invariants);

2. and that for a Kahler surface $S$ the Seiberg --- Witten invariant
for the canonical class $K_S$ equals to $\pm 1$,

(see [10]).

Based on this facts Witten reproved that a priori  Kahler surface $S$ isn't
diffeomorphic to its topomodel.

As in the Donaldson theory the invariants are defined if $b_2^+$ of based manifold bigger then 1 (for the Kahler case it just corresponds to $p_g (S) >0$).

But really the Seiberg --- Witten invariants use to be slightly more usefull 
then the Donaldson one. We have in mind the following important result due to
Taubes

\proclaim{Theorem (Taubes, [6])} Let $X, \om$ be  symplectic manifold with
$b^+_2(X) >1$. Then the Seiberg --- Witten invariant of the canonical class
$K_{\om} \in H^2(X, \Z)$, associated to given symplectic form, equals to
$\pm1$.
\endproclaim

In other words, for symplectic manifold the canonical class is a basic class.

But the inverse implication isn't true: in [3] the counterexample was found.
The construction is as follows. One can start with a symplectic manifold
$X$ with $b_2^+(X) >1$ and consider the connected sum of the following form
$$
Y = X \sharp N,
$$
where $b_1(N) = b_2^+ (N) = 0$ but $\pi_1(N) \neq \{1\}$. Then $Y$ doesn't admit
symplectic structure but has nontrivial invartiant (see [3]).

The idea to find a criterion for nontrivialty of the invariants was quite reasonable, and the next step in this direction was done in [9].
In this paper we've extended the Taubes technique and modulo some remarks, represented in the present paper, the following result
was established
\proclaim{Theorem} Let $X$ be pseudo symplectic manifold with $b_2^+(X) >1$, so there exists a hermitian triple $(g, J, \om)$ on $X$ with nontrivial
image in $H^2(X, \R)$ by the canonical map $\tau$. Then the invariant for
$K_J \in H^2(X, \Z)$ (which is, of course, the canonical class of $J$) is
nontrivial.
\endproclaim

Let us recall the definitions of pseudo symplectic manifolds  and the canonical map $\tau$. For a smooth 4 - manifold $X$ one can consider the space of all hermitian triples
$$
\sM_X = \{(g, J, \om)\},
$$
where the first element is a riemannian metric on $X$, the second is an almost complex structure, compatible with the first element, and the third is
the corresponding 2- form (which is called an almost Kahler form) such that
$$
\om (u, v) = g(u, Jv).
$$ 
It's well known that this form $\om$ is

1. self dual with respect to the conformal class $*_g$ and the orientation,
choosen by $J$;

2. nondegenerated everywhere;

3. has the type (1,1) with respect to $J$.

We have to emphercise
that this $\om$ hasn't to be closed a priori. Because of this it is just
{\it  almost Kahler}.

The construction is pure local and it's well known from the linear algebra that in a hermitian triple every element can be reconstructed from two others.

Further, there exists the canonical map 
$$
\tau: \sM_X \to K^+ \subset H^2(X, \R),
$$
where $K^+$ is the inner part of  isotropic cone (so 
consists of  2 - cohomology
classes with positive squares), defined as follows. For a triple
$(g, J, \om)$ let us take the corresponding Hodge star operator 
$$
*_g: \Om^i_X \to \Om^{4-i}_X,
$$
over $X$ defined by the conformal class of our riemannian metric $g$ and the orientation given by $J$.
Then by the famous Hodge theorem ([2])
 every form can be decomposed into three
parts --- harmonic part, exact and co- exact 
$$
\om = \om_H + d \rho_1 + d^* \rho_2,
\tag 0.1
$$
where $\rho_1 \in \Om^1_X, \rho_2 \in \Om^3_X$.

Let us recall that $d^*$- operator on $\Om^2$ is equal to  the following double combination
$$
d^*: \Om^2 \to \Om^1 \quad | \quad d^* = *_g d *_g
$$
of the ordinary $d$ and two Hodge operators (in different dimensions).

But in our case $\om$ is self dual, so one can apply $*_g$ to both sides of
(o.1) and get 
$$
*_g \om = \om = *_g \om_H + *_g d \rho_1 + d (*_g \rho_2).
$$

From the uniqueness of the decomposition one imidiatelly gets that
$$
\rho_1 = *_g \rho_2
$$
and 
$$
*_g \om_H = \om_H,
$$
hence one can rewrite the previous decomposition as follows
$$
\om = \om_H + d^+\rho_0, \quad \rho_0 \in \Om^1_X \quad {\text where}\quad
d^+: \Om^1_X \to \Om^+_X, 
\tag 0.2
$$
and establish that the harmonic part $\om_H$ is self dual too.

\proclaim{Definition} The image of the canonical map $\tau$ for  given
hermitian triple $(g, J, \om)$ equals to the 2- cohomology class which is
represented by the harmonic form $\om_H$ defined in (0.1), (0.2)
$$
\tau(g, J, \om) = [\om_H] \in H^2(X, \R).
\tag 0.3
$$
As we've seen above $\om_H$ is self dual so it has positive square
$$
\int_X \om_H \wedge \om_H = \int_X \vert \om_H \vert_g^2 dm_g > 0
$$
unless the case when $\om_H$ is vanishing everywhere. Hence the image
$Im \tau$ lies in $K^+$.
\endproclaim 

\proclaim{Definition} One can called based manifold $X$  pseudo simplectic if there exists 
a hermitian triple with nontrivial image in $K^+$. \endproclaim

And it was proved that the canonical class of {\it pseudo symplectic manifold}
has the invariant equals to $\pm1$ as well as in the original symplectic case.

\subheading{Example} First of all it's clear that every symplectic manifold is pseudo symplectic. For a given symplectic manifold $X$ endowed with a symplectic form $\om$ one gets that for an apropriate riemannian metric $\om$ is harmonic itself so $\om_H = \om$ and hence $\tau(g, J, \om) = [\om] \ne 0$. But we'd like to recall an example of pseudo symplectic
manifold which doesn't admit a symplectic structure. Namely, if one take
the counterexample from [3] it isn't too hard to establish that
this connected sum $Y = X \sharp N$ (see [3], [9]) admits 
a hermitian triple with nontrivial image in $K^+$. Such triple on whole $Y$
can be constructed as an extension of the symplectic triple on the original
symplectic manifold $X$ (we'll discuss all details in the next section).

After all above is understood there is one absolutely reasonable quastion. It's well known that in our setup the Seiberg --- Witten invariant depends only on the choice of $Spin^{\C}$ - structure so on the corresponding class $c \in H^2(X, \Z)$ and doesn't depend on the choice of riemannian metric and almost complex structure. But from the first viewpoint two hermitian triples $(g_i, J_i, \om_i), i = 1, 2$ with the same canonical class $K= K_{J_i}$ can have
rather different images in $K^+$. The first can be nontrivial whether the second
is trivial.

In this paper we proof the following 
\proclaim{Main Theorem} The nontriviality condition for image of $\tau$ is stable with respect to contineous deformations. It means that if $\tau (g_0, J_0, \om_0) \ne [0] \in H^2(X, \R)$ the same is true for every $(g, J, \om)$
with $K_J = K_{J_0}$.\endproclaim

On the other hand we can reformulate the statement with respect to
the action of the group
$$\sG = Aut(TX)
$$
of all smooth fiberwise automorphisms of the tangent bundle over $X$.
Namely there exists unique discrete invariant of triples --- the corresponding
canonical class $K_J \in H^2(X, \Z)$. And it's more or less clear that
if two triples $(g_i, J_i, \om_i), i= 1, 2,$ have the same
canonical class then there exists an element $u \in \sG$ which conjugates
the first triple to the second. So as well one can say that the nontriviality
of the image of $\tau$ {\it is stable with respect to the $\sG$- action}.

\head \S 1. Pseudo symplectic manifolds in the framework
of Seiberg --- Witten theory
\endhead

In this section we'll sketch the proof of the Taubes Theorem
mentioned above and show in which point it uses
symplectness of based manifolds.

First off all in absolutely general almost complex situation
for a hermitian triple $(g, J, \om)$ one has the following objects:

1. the decomposition on self dual and anti self dual forms
represents as
$$
\aligned \La^+ = &(\La^{2,0} \oplus \La^{0,2})_{\R} \oplus \R <\om>, \\ \La^- & = (\R <\om>)^{\perp} \subset \La^{1,1};
\endaligned
\tag 1.1
$$

2. spinor bundles correspond to $Spin^{\C}$ - structure
$c_0 = - K_{\om} \in H^2(X, \Z)$ have the form
$$
\aligned
W^+ = & \La^{0,0} \oplus \La^{0,2} = I \oplus K^{-1},\\
W^- & = \La^{0,1}\\
\endaligned
\tag 1.2
$$
and our self dual form $\om$ (as the section of
$\La^+ \cong ad W^+$) acts on $W^+$ 
as diagonal operator with eigenvalues $2 \imath$
and $- 2 \imath$ on the direct summands $K^{-1}$ and
$I$ respectively;

3. the gauge group is
$$
\sG = Aut_h (det W^+) = Aut_h (K^{-1});
\tag 1.3
$$

4. the configuration space for the Seiberg --- Witten system
is represented by space of triples $(a, \al, \be)$
where $a$ is a hermitian connection on $K^{-1}$,
$\al$ is a section of $I$ so is a complex valued function
and $\be$ is a section of $K^{-1}$;

5. the Seiberg --- Witten system in this setup reads as
follows
$$
\aligned
& D_a (\al \oplus \be) = \overline{\partial}_a \al + \overline{\partial}_a^* \be = 0 \\
& \imath F^+_a = (\vert \be \vert^2 - \vert \al \vert^2) \om
+ \al \overline{\be} + \overline{\al} \be; \\
\endaligned
\tag 1.4
$$

and let us recall that in this case the moduli space of solutions is zero dimensional.
(For all details see [4], [6], [9], [10].)

\subheading{ The first step}
In any case there exists unique up to gauge transformations
hermitian connection $a_0 \in \sA_h (K^{-1})$ such that
the projection of the corresponding connection (on whole 
$W^+ = I \oplus K^{-1}$ ) to the first direct summand
is equal to ordinary $d$
$$
\nabla_{a_0} |_I = d \quad {\text on } \quad I.
\tag 1.5
$$
So if one take spinor filed of the form $(s \oplus
0) \in I \oplus K^{-1}$ where $s$ a constant function then 
$$
\nabla_{a_0} (s \oplus 0) = b \in \Ga(K^{-1} \otimes T^*X)
\tag 1.6
$$
and $b$ essentually is {\it the torsion of our almost complex structure} $J$. So in the case of integrable complex structure one has $b = 0$ and hence $D_{a_0}(s \oplus 0)$
automatically equals to zero. But in the symplectice case
Taubes observes that despite of nontriviality of
$\nabla_{a_0} (s \oplus 0) = b$ the corresponding Dirac operator $D_{a_0}$ vanishes on $s \oplus 0$ if and only
if the form $\om$ is closed (so is a symplectic form).

Taubes arguments on this step are  extremly usefull and
deep exercises in the theory. Namely he acts by
the self dual form $\om$ on the spinor field $s\oplus 0$
$$\om (s \oplus 0) = -2 \imath \cdot (s \oplus 0)$$
and then applies the corresponding coivariant deravitive
to the both sides of the previous equality
$$
\nabla_{a_0} (\om (s \oplus 0)) = - 2 \imath b.
$$
But on the right side we have by the Liebnitz rule
$$
\nabla_{a_0} (\om (s \oplus 0)) = C \cdot d^* \om \otimes (s \oplus 0) + \om (b) = C \cdot d^* \om \otimes (s \oplus 0) + 2 \imath b
$$
(where $C$ is an integer number) hence   using Clifford multiplication one gets
$$
D_{a_0} (s \oplus 0) =   \frac{\imath C \cdot s}{4}
 (d^* \om)^{0,1}
\tag 1.7
$$
where the last term $(d^* \om)^{0,1}$ can be equal to zero
if and only if $d^* \om$ is equal to zero so if and only if
our form $\om$ is closed and hence symplectic.

\subheading{The second step} After it was established
that the spinor field $s \oplus 0$ is harmonic with respect
to the Dirac operator $D_{a_0}$ Taubes considers some special perturbation of the original Seiberg --- Witten system
namely
$$
\aligned
& D_{a_0} (\al \oplus \be) = 0 \\
& \imath F^+_a =  \imath F^+_{a_0} + (\vert \be \vert^2 -
\vert \al \vert^2 + 1) \om + \al \overline{\be}
+ \overline{\al} \be. \\ \endaligned
\tag 1.4'
$$
It's clear that the triple $(a = a_0, \al = 1, \be = 0)$
is a solution of this perturbed system. But the point is
that one can derive the Seiberg --- Witten invariant from this perturbed system as well as from the original one.
It's so since {\it there are no reducible solutions}
for the perturbed system. Really for such solution (with trivial spinor part) one has 
 $$
\imath F^+_a = \imath F^+_{a_0} + \om
$$
but these two curvature forms are cohomologically the same so
$$
\om = d^+ \rho
$$
where $\rho$ is a real 1- form. But so as
we are now in the symplectic cathegory
it's impossibly hence there are no reducible solutions.

\subheading{The third step} Then Taubes imposed an additional term to the right side of the second equation
$$
\imath F^+_a = \imath F^+_{a_0} +
(\vert \be \vert^2 - \vert \al \vert^2 +1) \om
+ \al \overline{\be} + \overline{\al} \be + \frac{4r}{1+r \vert \al \vert^2}
(\bar \al <b, \nabla_{\xi} \al> + \al \overline{<b, \nabla_{\xi} \al>}).
\tag 1.8
$$
Here $<., .>$ is the $\C$- bilinear extension to $T^*_{\C}X$ of the metric inner product, $b$ is the above- defined section of $K^{-1} \otimes T^*X$, and
$\nabla_{\xi} = d + \imath \xi$ is a covariant derivative of complex- valued functions on $X$, where the pure imaginary form $\xi$ is the difference $\xi = a - a_0$ of two Hermitian connections from $\sA_h(K^{-1})$ (see [6], [4]).
Then it was  established that for sufficienly large $r >> 0$
there exists unique up to gauge transformations solution
with multiplicity one for this perturbed system and hence the invariant equals
to $\pm 1$ (see [6], [4]).

We recalled more or less rigouresly just the first and the second steps because only on these steps Taubes uses
{\it closedness} of the self dual 2- form $\om$.

So to generalize this imroptant result one has to avoid
in these steps contradictions and to impose a condition on
hermitian triples.

\subheading{The first step}
We've seen that in absolutely general situation (in the framework of the almost complex geometry)
$$
D_{a_0} (1 \oplus 0) = \frac{\imath C}{4}
 (d^* \om)^{0,1} \in \Om_X^{0,1}
= \Ga (W^-)
$$
supposing that $s = 1$. Then one can deform the connection
$a_0$ such that the corresponding covariant derivative
has the form
$$
\nabla_{a_1} = \nabla_{a_0} +  \pmatrix
\frac{C}{4 \imath} & 0 \\ 0 &  \frac{C}{4 \imath} \\
\endpmatrix  \otimes d^* \om \in \sA_h (W^+).
\tag 1.9
$$
It corresponds to such hermitian connection $a_1 \in \sA_h (det W^+)$ that
$$
a_1 = a_0 + \frac{C}{4 \imath} d^* \om
\tag 1.10
$$
and it's clear so as the last term is a pure imaginary
1- form that such deformation is well defined.

Then we can use this new connection $a_1$ instead of the original one in absolutely general case to establish
that
$$
D_{a_1} (1 \oplus 0) = 0
$$
always (direct substitution of (1.9) to (1.7)).

\subheading{The second step}
So on this level one can use a familiar perturbation
of the original system which differs from the Taubes perturbation just in the second equation
$$
\imath F^+_a = \imath F^+_{a_1} + (\vert \be \vert^2
- \vert \al \vert^2 + 1) \om +
\al \overline{\be} + \overline{\al} \be.
\tag 1.8'
$$
Therefore again we have "primitive" solution which is
now $(a = a_1, \al = 1, \be = 0)$
and one has to establish non existance of reducible solutions
for our perturbed system.

At this point our definition of pseudo symplectic triples
has to be imposed. Namely let us suppose
that a reducible solution exists for our choosen triple
$(g, J, \om)$. Then we get that
$$
\imath F^+_a - \imath F^+_{a_1} = \om
$$
but again two curvature forms in the left side are cohomologically the same so in this case it should be
a real 1- form $\rho$ such that
$$
\om = d^+ \rho.
\tag 1.9
$$
But such $\rho$ can exist if and only if
the corresponding harmonic part $\om_H$ in the Hodge decomposition (0.1)
$$
\om = \om_H + d \rho_1 + d^* \rho_2
= \om_H + d^+ \rho_0
$$
is trivial. So if our hermitian triple $(g, J, \om)$ has nontrivial image $\tau (g, J, \om) \ne [0]$
in $H^2(X, \R)$ (see the Introduction)
there are no reducible solutions for our perturbed system. Then consider
the following version of the SW - system, combined form (1.8) and (1.8'): 
$$
\aligned
D_a(\al + \be) = 0, \quad \imath F^+_a = \imath F^+_{a_1}
+(\vert \be \vert^2 - \vert \al\vert^2 + 1) \om + \\
\al \bar \be + \bar \al \be +
\frac{4 r}{1 + r \vert \al \vert^2}(\bar \al <b, \nabla_{\xi} \al>
+ \al + \overline{<b, \nabla_{\xi} \al>}),\\
\endaligned
\tag 1.8"
$$
where all the notation are explained above in (1.8). So it differs from
the Taubes $r>0$ perturbation (see [4], [6]) just by the self- dual
form $\imath F_{a_1}^+ - \imath F_{a_0}^+$ in the second equation.
Now we claim that our basic solution $(a_1, 1 \oplus 0)$ satisfies
the both equations from (1.8"). Really, to check this one could 
test only the fifth term of the second equation right side. 
So it's sufficient to check that 
$$
<b, a_1 - a_0> = 0.
$$
But by the construction above the difference $a_1 - a_0$ is proportional
to $d^* \om$. And moreover, the same construction from the Step One gives us that our $b$ is related to this real 1 - form. Namely 
$$
b \in \Ga (\La^{0,2} \otimes T^*X) \mapsto b' \in \Ga(\La^{1,2})
\mapsto \rho_b \in \Ga(T^*_{\C}),
$$
where the last map is defined by the operator $\La$, adjoint to the wedge
product by $\om$. By the construction $\rho_b$ is proportional to
$(d^* \om)^{0,1}$ (see the Step One), so using local computation one gets
that 
$$
<b, \rho_b> \in \Ga(\La^{0,2})
$$
{\it is trivial} for every section of $K^{-1} \otimes T^*X$ and hence it is trivial for our $b$.

Thus our $(a_1, 1 \oplus 0)$ is a solution of (1.8") for every $r\geq 0$.
Then one can repeat all Taubes calculation to establish that for $r>>0$
there exits only one solution of the system (1.8") and this gives us
that the Seiberg --- Witten invariant is non trivial for $Spin^{\C}$ - 
structure $- K_J$. We do not discuss the question of the multiplicity
of this solution, because rather then in the symplectic case one cann't claim
that it  equals to one.

Let us continue with the example of pseudo symplectic but non symplectic manifolds mentioned in the Introduction.
So $X$ is a symplectic manifold with $b^+_2 > 1$ and symplectic form $\om$. Consider the following
connected sum
$$
Y = X \sharp N
$$
where $N$ is a manifold with $b_1(N) = b_2^+(N) = 0$.
It was shown in [3] that such $Y$ has the same invariants
as the original $X$ but doesn't admit symplectic structure
if the fundamental group $\pi_1(N)$ admits nontrivial finite quotient.
To establish that there are no symplectic structures on $Y$ it's sufficiently to consider the universal covering
of $Y$ which is a smooth manifold of the following form
$$
\wave{Y} = X \sharp X \sharp ... \sharp X \sharp N'
$$
and here we have $d$ copies of $X$ where $d$ is the corresponding order of the quotient. Since any symplectic structure on $Y$ induces the corresponding symplectic structure on $\wave{Y}$ combinig this with the fact
that every connected sum where summands have positive 
$b^+_2$ has trivial Seiberg --- Witten invariants 
one gets that $Y$ doesn't admit any symplectic structure.

But for our aims such examples are quite appropriate. Let us formulate the result as follows
\proclaim{Lemma} Let $X$ be symplectic  4- manifold endowed with symplectic form $\om_0$. Then for every 4- manifold $N$ 
such that $b_1(N) = b_2^+(N) = 0$ the connected sum
$$
Y = X \sharp N
$$
is a pseudo symplectic manifold.
\endproclaim

To proof this Lemma one has to construct a pseudo
symplectic structure on $Y$ so to find a hermitian triple
$(g, J, \om)$ with nontrivial image in $H^2(Y, \R)$.

Let us start with the given symplectic structure
on $X$ represented by a triple $(g_0, J_0, \om_0)$
where the third element is our given symplectic form
and we choose some compatible riemannian metric
$g_0$ to form such triple.  Then let us
extend abritrary this original triple to the whole
$Y$, getting some triple $(g_1, J_1, \om_1) \in
\sM_Y$. First of all the corresponding canonical class
$$K_{J_1} \in H^2(Y, \Z) = H^2(X, \Z) \oplus
H^2(N, \Z)$$ has the following form
$$
K_{J_1} = K_{J_0} +   \sum_{1}^{b_2(N)} c_i
\tag 1.10
$$
where $c_i$ is the basis in $H^2(N, \Z)$ in which
the intersection form $Q_N$ has standard form
$(-1,..., -1)$.  Obviously we can extend
all the elements if such choise is done.
Really,  let us realize our almost complex structure
$J_0$ as a section of the corresponding spinor bundle
$W^+_X \to X$ with the Chern classes
$c_1(W^+_X) = K_{J_0}, c_2 (W^+_X) = 0$.
Then it's clear that over $N$ the sum $  \sum_{1}^{b_2(N)}
c_i$ is a $Spin^{\C}$- structure (because it is a characteristic element of $Q_N$) and one can take
the corresponding spinor bunle $W^+$ over $N$ with the
Chern classes $c_1 (W^+_N) =   \sum_{1}^{b_2(N)} c_i,
c_2(W^+_N) = 1$ induced by some appropriate riemannian metric
over $N$. Then over the whole $Y$ we get the spinor bundle
$W^+_Y$ with the Chern classes $c_1(W^+_Y) = K_{J_1},
c_2(W^+_Y) = 0$ (it's easy to check it directly). So one can 
extend the choosen section of $W^+_X$ corresponded to our given
almost complex structure and then get a  nonvanishing section of $W^+_Y$ which would induce an almost complex structure on the whole $Y$.

So we claim now to deduce that the triple $(g_1, J_1, \om_1)$
has nontrivial image $\tau(g_1, J_1, \om_1)$ in $H^2(Y, \R)$.
For this let us take the harmonic with respect to $*_{g_1}$
2- form $\om_H$ such that
$$
[\om_H] = [\om_0] \in H^2(Y, \R).
$$
Then consider the following integral
$$
\int_Y \om_1 \wedge \om_H.
\tag 1.11
$$
Since  both 2- forms are self dual with respect to
$*_{g_1}$ this expression equals to {\it the inner product}
so if it is nontrivial then the projection of $\om_1$
to the ray $\R <\om_H>$ in the harmonic space is nontrivial
and hence $\tau (g_1, J_1, \om_1)$ has to be nontrivial too.

As usual for connected sums let's divide
the integral (1.11) into three parts
$$
\int_Y \om_1 \wedge \om_H =
\int_{X \backslash B_1} \om_1 \wedge \om_H +
\int_{{\text Neck}} \om_1 \wedge \om_H +
\int_{N \backslash B_2} \om_1 \wedge \om_H
$$
where $B_i$ are  small balls using for the glueing
procedure.

Then since over punctured $X \backslash B_1$
these two 2- forms are very closed to the original symplectic form
we get that 
$$
\int_{X \backslash B_1} = 2 Vol X + o(r),
$$
where $r$ is the radius of the balls. So as
the Neck is conformal flat one gets
$$
\int_{{\text Neck}} \om_1 \wedge \om_H =
o(r),
$$
and since there are no self dual harmonic form
over $N$ the same expression can be obtained for the third summand
$$
\int{N \backslash B_2} \om_1 \wedge \om_H
= o(r).
$$
Taking all together we see that
$$
\int_Y \om_1 \wedge \om_H =
2 Vol X + o(r)
\tag 1.12
$$
so shrinking the Neck (standard trick, see f.e. [1])
we get that this integral is nontrivial and positive.

The Lemma proved above gives  us 
a number of examples of pseudo symplectic but non symplectic
4- manifolds.

This story  represents the motivation to introduce the definition of
pseudo symplectic manifolds. Now  we'd like to repeat
the natural conjecture which was stated above:
we know that the invariants doesn't depend
of the choice of riemannian metric and depends only on the corresponding class  in $H^2(X, \Z)$. Is the same true for  nontriviality of the image of our canonical map $\tau$?

The answer (included in our Main Theorem) is positive and
in the following sections we prove this statement.

\head  \S 2. The structure of the space $\sM_X$ for a given $X$
\endhead

In this section we'll describe the structure of the hermitian triple space for a given smooth 4 - manifold 
$X$. 

First of all let's consider the discrete invariants, which divides all the $\sM_X$
into discrete set of subspaces.

The most generic division of $\sM_X$ corresponds to the orientation. For an orientable $X$ there are two choices of orientation. And if one fix an almost complex structure then the corresponding orientation is fixed. So  one can divide $\sM_X$ into two pieces
$$
\sM_X = \sM_X^+ \cup \sM_X^-,
\tag 2.1
$$
where triples from the first subspace are compatible with one orientation
and triples from the second one are compatible with the other.
 There are not any preferences in the choice of the orientations but let us assign $+$ and $-$ arbitrary
 to
distinguish the two orientations.

The next (and the last) discrete data are the labelling of $\sM_X$ by the corresponding canonical classes. Namely let's consider the following obvious map
$$
Can: \sM_X \to K_X \subset H^2(X, \Z)
\tag 2.2
$$
defining by
$$
Can(g, J, \om) = K_J \in H^2(X, \Z)
$$
where $K_J$ is the canonical class of $J$. There are two aproaches to define
this canonical class. First of all, the tangent bundle $TX$ together with
our choosen almost complex structure $J$ can be regarded as rank two complex bundle $T^{1,0}_{\C} X$, so one can take the first Chern class of this bundle
and
$$
K_J = - c_1 (T^{1,0}_{\C}X).
$$
The second approach is to consider the direct sum
$$
\La^2 (T^*X) = \La^+ \oplus \La^-
$$
where $\La^{\pm}$ are real rank- 3
bundles on self dual (respectively anti self dual) 2- forms with respect to
the conformal class $*_g$ and the orientation defined by $J$. It's well known
that there is trivial real subbundle $\R <\om> \subset \La^+$ and in the presistance of the riemannina metric $g$ one gets the ortogonal to $\R<\om>$
subbundle $(\R<\om>)^{\perp} \subset \La^+$. Moreover, the vector multiplication by
our form $\om$ defines a complex structure on $(\R<\om>)^{\perp}$ and one can take
the first Chern class of this line complex bundle and it is our canonical
class again.

It's easy to decribe the image of $Can$ in $H^2(X, \Z)$. Let us fix an orientation
on $X$, then for any $K_0 \in Im Can \subset H^2(X, \Z)$ one has
$$
K_0 = w_2 (X)(mod 2), \quad K_0^2 = 2 \chi + 3 \si,
\tag 2.3
$$
where $w_2(X)$ is the second Stiefel - Whitney class, $\chi$ is the euler
characteristic (and these elements are independent of the choice of orientation)
and $\si$ is the signature, corresponding to the choosen orientation. 

Further, for any $K_0 \in Im Can$ one has that $-K_0 \in Im Can K$. Namely for any triple
$(g, J, \om)$ with $K_0$ we have $(g, -J, -\om)$ with $-K_0$. So
the space $\sM_X$ has {\it standard real structure} $\Theta_X$ without real points. For any point $(g, J, \om) \in \sM_X$ we have
the conjugated point 
$$\Theta_X (g, J, \om) = (g, -J, -\om) \in \sM_X \tag 2.4$$ and obviuosly
our real structure $\Theta_X$ preserves both $\sM_X^{\pm}$ and induces
just the multiplication by $-1$ on $K_X \subset H^2(X, \Z)$.  

Hence our space $\sM_X$ is labelled by the discrete image of the map $Can$ and we have
$$
\sM_X / \Theta_X = \bigcup_{i} \sM_X^{K_i},
$$
where $\sM_X^{K_i} = Can^{-1}(K_i)$.

To describe the structure of the space $\sM_X$ we
 consider a triple $(g, J, \om)$ with fixed canonical class
as the corresponding pair $g, J$. Then there are two natural projections:
to the first and  the second components
$$
\pi_1 (g, J) = g, \quad \pi_2 (g, J) = J.
\tag 2.5
$$
The images lie in the homogeneous space of all riemannian metrics $\Cal H$, compatible with the given smooth structure, and in the space of all almost complex structures with the fixed canonical class. Now we want to describe
the fibers of these two projections. But for this
at a moment we'll use not the space $\sM^{K_i}$ but the space of pairs
$(*_g, J)$ where $*_g$ is the conformal class of a riemannian metric. Let us
denote this space as $\wave{\sM^{K_i}}$, so we have the following natural fibration
$$
con: \sM_X^{K_i} \to \wave{\sM_X^{K_i}}
$$
which is a principal $e^{\R}$- bundle. 

To describe the  projections above it's very usefull to translate all constructions to the language of projectivization. In this setup
the conformal class $*_g$ of a metric $g$ corresponds to
the following object. Over each point $p \in X$ one has in 
$$
\C\proj^3 = \proj(T_p^{\C}X)
$$
non- degenerated quadric $Q_p \subset \C \proj^3$ together with standard
real stucture $$\Theta_p: \C\proj^3 \to \C\proj^3, \quad \Theta_p^2 = id$$
such that $Q_p$ is real but without real points. Then any almost complex structure $J$ compatible with this conformal class corresponds to a pair
of projective lines $l_p^1, l_p^2$  on the quadric. It's well known that there are
two families of projective lines on a nondegenerated quadric. And our
two lines lie in the same family because
$$
\Theta_p (l_p^1) = l_p^2
$$
and if $l_p^1, l_p^2$ will be in two different families 
then the intersection point would be real. Our nondegenerated quadric
$Q_p$ is the direct product of two projective lines 
$$
Q_p = \proj^+ \times \proj^-
$$
(since the orientation is choosen by the fixing of $K_i$ one can distinguish
these two projective lines by the signs). So by the definition $l_p^{i}$  is represented by  points in $\proj^+$ (for details see [10]).

In this language the description of $\pi_i$ is quite obviuos.

\subheading{ The first projection}
Now we're ready to describe a fiber of
$$
\pi_1 : \sM^{K_i} \to \Cal H,
$$
where $\Cal H$ is the homogeneous space of all riemannian metrics compatible with the given smooth structure.

Let us fix a riemannian metric $g$ and a canonical class $K_i$. Then one gets
the following lifting of the corresponding projective bundle $\proj^+$.  Since our
$X$ is orientable (and the corresponding orientation is fixed by our canonical class) the  lifting is defined by the choice of the first Chern class (which is called the choice of $Spin^{\C}$ - structure). So if we choose $K_i$ as the first Chern class we immediately get the corresponding vector bundle $W^+$ which is called {\it spinor bundle}. The topological type of the bundle is defined by
our choosen $K_i$
$$
c_1(W^+) = K_i, \quad c_2(W^+) = 0.
$$
Then one can consider the space of smooth everywhere nonvanishing sections: 
$$
\Ga_*(W^+) = \{ \phi \in \Ga(W^+) \quad | \quad (\phi_0)|_0 = \emptyset\}.
\tag 2.6
$$
Then any section $\phi \in \Ga_*(W+)$ defines an almost complex structure as follows. There is natural pairing $\phi \to (\phi \otimes \overline{\phi})_0$
and in the composition with well- known isomorphism
$$
ad W^+ = \La^+
$$
one gets the corresponding self dual 2- form $\om_{\phi}$ which is nondegenerated everywhere. So from $g$ and $\om_{\phi}$ one can reconstruct
the second element in hermitian triple and it is our almost complex structure
(see, f.e., [8]).

But it's clear that two nonvanishing sections define the same complex structure
if and only if they are {\it gauge equivalent} with respect to the gauge group $$\sG_0= Aut_h (det W^+)$$ of fiberwise transformations of the determinant line bundle $det W^+$.
So a fiber of $\pi_1$ above is naturally isomorphic to
$$
\pi_1^{-1} = \Ga_*(W^+) / \sG_0.
\tag 2.7
$$

Equivalently one can describe the first projection in terms of the space
of nonvanishing self dual 2- forms. Again, since  we have choosen orientation and
riemannian metric $g$, we can decompose the space $\Om^2_X$ of all smooth 2- forms into self dual and anti self dual parts
$$
\Om^2_X = \Om^+ \oplus \Om^-.
$$
Then let us derive the subset $\Om^+_* \subset \Om^+$
of nonvanishing everywhere self dual forms. Let us recall that for
self dual 2- forms there are just two possibilities: to have rank 0
or rank 4 (full rank). So the subset $\Om^+_*$ is defined by just one condition:
to be nonvanishing. Then the subset $\Om^+_*$ consists of a number of
disconnected components, labelled as above by the corresponding canonical classes. In this sense the present picture is more universal then the picture above (with the spinor bundle) because now  all the components of $\pi_1^{-1}(g)$ are defined on the  whole space $\sM_X$. So the corresponding connected component of
$\Om^+_*$ is the fiber of $\pi^{-1} \in \sM_X^{K_i}$.

On the other hand we can realize the situation globally using the following
global construction. For our fixed riemannian metric $g$ we consider
pointwise all self dual 2- forms which  correspond to almost Kahler forms.
Namely over each point $p$ we have a point $\om_p$ in the fiber $\La^+_p$.
But one has to unify all these forms so that each $\om_p$ has the norm
equals to $\sqrt{2}$. Hence we have 2- sphere $S^3_p \subset \La^+_p$ consists
of locally compatible forms. After globalization one gets the following fibration
$$
s: \Cal S \to X
\tag 2.8
$$
with a fiber $s^{-1}(p) = S^2$. It's clear that

a) $\Cal S$ is smooth compact orientable 6- dimensional  manifold;

b) the topological type of $\Cal S$ is fixed by the topological type of $X$;

c) every global everywhere nonvanishing  self dual 2- form $\om$ defines
smooth inclusion
$$
i_{\om} X \hookrightarrow \Cal S
$$
such that the image $i_{\om}(X)$ intersects all the fibers transversally;

d) the 4- homology classes of the images $[i_{\om_i}] \in H_4(\Cal S, \Z)$ are in
one - to - one correspondence with the canonical classes defined by $\om_i$.

Moreover since every 2- sphere fiber is endowed with the natural riemannian metric 
one gets the corresponding riemannian metric $G_g$ on whole $\Cal S$.
So:

e) the projection of $G|_{i_{\om}(X)}$ is isomorphic to $g$ on $X$.

And we can define an almost complex structure $\Cal J_g$ on $\Cal S$ as follows.
Over each point $t \in \Cal S$ there is ortogonal decomposition of the tangent space induced by our metric $G_g$
$$
T \Cal S_t = TS^2 \oplus TX
$$
and we  take the standard complex structure induced by the vector multiplication by $\om$ on the first summand and  the  tautological  almost complex structure defined by $g$ and $\om$ on the second summand.

Really the space $\Cal S$ is well known in differential geometry as twistor space but originally it was defined for conformal classes instead of
riemannian metrics themselfs. But the correspondence is absolutely clear
and will be used further in this paper.

\subheading{ The second projection}
Now let's consider the second projection
$$
\pi_2 :\sM_X^{K_i} \to \Cal K_{K_i}
$$
to the space of almost complex structure with the same canonical classes.

For this we'd like to present the {procedure to reconstruct 
riemannian metrics} compatible with a given almost complex structure from
some additional data. 

First of all let us recall that our standard real structure $\Theta_p$
on $\C\proj^3_p$ induces  {quaternionic real structures} on our
projective lines $\proj^{\pm}_p$ which parametrize the two families on the quadric $Q_p$.

Suppose that we already fixed an almost complex structure $J$ and got over each point 
$p \in X$ two projective lines $l^i_p$. 
The same picture exists over each point and one can see that $J$ defines
an inclusion
$$
\rho_J:\proj^- \hookrightarrow \C\proj^3
$$
of projective bundle (global!) $\proj^- \to X$ to our projective bundle
$\C\proj^3 \to X$ so the image is a projective subbundle of $\C\proj^3$.
Over each point 
$$\rho_J(\proj^-)|_p = l^1_p. \tag 2.9$$

Then we can reconstruct a quadric
$Q_p$ from the corresponding quaternionic real structure $\theta_p$ on
$\proj^-_p$ by the following procedure. Together with the identification
(2.9) above one gets the corresponding  quaternionic real structure 
$$ 
\theta_p' = (\rho_J)_*(\theta_p)$$ on $l^1_p$. We define for each point $s^1 \in l^1_p$
a point $s^2 \in l^2_p$ as
$$
\Theta_p (\theta_p'(s^1)) = s^2.
$$
Then if we take all projective lines of the form
$$
<s^1, s^2>
$$
we'll get a nondegenerated quadric which is of course our $Q_p$.
 An extra datum for the reconstruction procedure is the following: over each point we have a global pintwise quaternionic structure on $\proj^-$. So the space of
all conformal classes (we have to emphercise thet we get just a conformal class) compatible with a given almost complex structure $J$ can be regarded as
the space of all global sections  of the corresponding principal $PGL(2, \C)$- bundle.
 After local consideration we again insure that $Aut(TX)$ acts transitively
on $\sM_X^{K_i}$. But  it's clear that
for every global $\theta$ there is nontrivial stibilizer in $PGL(2, \C)$.

Let us note that really the picture of the reconstruction is more {\it universal} in the following sense. For a fixed
quaternionic real structure on the projective bundle 
$\proj^-$  one can get different metrics   using different almost complex structures
with different (ad hoc) canonical classes. So  the unique
fixed additional datum  corresponds to the subset in the whole $\sM_X$.

The existance of the stabilizers  gives us the following conclusion. The space $\sM_X$ has no  structure of a principal bundle but is
a homogeneous space. This is easy to see from the consideration of any of our two projections.
So let's now consider the fibration discribed above
$$
Can: \sM_X \to K_X.
$$
Then we have the action of $Aut(TX)$ on each fiber. But for every point
$(g, J, \om) \in \sM_X$
there is nontrivial stabilizer which is the space of sections of a principal $U(2)$- bundle. We'd like to postpone more precise description to a future
since our goal now is just to study the pseudo symplecticity condition.

\head \S 3. Inclusions into twistor space
\endhead

We can discribe a fiber of the first projection $\pi_1$  as it was outlined
above using smooth inclusions and turning to an analogy with complete linear systems in the algebraic geometry.

Let us fix a conformal class $*_g$ and consider the corresponding fiber
of $\pi_1$. Then over each point we have $\proj^+_p$ projective lines
so that our quadric $Q_p$ is the direct product of two   projective lines
parametrizing two families of projective lines on $Q_p$. Let us globalize
this picture over whole $X$ to get a projective bundle $\proj^+ \to X$.
The total space of this bundle is {\it twistor space} of the conformal class
which we'll denote as $Y$. The topological type of this $Y$ is fixed
by the topological type of $X$ and moreover it has a smooth structure
and an almost complex structure $\Cal J_g$ both defined by the smooth structure on $X$ and our conformal class. We have 
$$
\pi: Y \to X
\tag 3.1
$$
such that $\pi^{-1}(p) = \C\proj^1$
which is our projective line above. So for any compatible
almost complex structure $J$ on $X$ one has an inclusion
$$
i_J: X \to Y
\tag 3.2
$$
such that

1. $i_J$ is smooth inclusion;

2. the projection $d \pi (\Cal J_g  |_{i_J(X)}) = J$ by the definition;

3. the corresponding 4- homological class $[i_J(X)] \in H^4(Y, \Z)$ depends only on the corresponding canonical class $K_J$.

Here we need some computation to find what it the homological class
$[i_J(X)]$ in $H_4(Y, \Z)$. For this one can note that {\it topologically}
our twistor space is isomorphic to the projectivization $\proj(W^+)$ of the spinor bundle $W^+$. Then we use the standard knoledgements about 
projectivizations of complex bunldes over $X$. Namely if $E$ is a complex rank 2
bundle over $X$ then the projectivization $\proj(W^+)$ has the following
cohomology ring
$$
H^2(X, \Z) \oplus \Z [H] \backslash (c_2(E) + c_1(E) \cdot H + H^2).
$$
Of course, this generator $H$ is defined modulo twisting by complex line bundles. So in our case with $E = W^+$ we choose $W^+ = I \oplus K^{-1}_J$
and thus we get that in this case the homology class $[i_J(X)]$ is Poincare dual
to the class $[H]$. Moreover, for any other $J_1$ we have
$$
[i_{J_1}(X)]^* = [H] + \half(\pi^*(K_{J_1} - K_J)).
$$
For example, for the almost complex structure $-J = J^{-1}$
with canonical class $-K_J$ we get
$$
[i_{-J}(X)] = [H] - \pi^*(K_J),
$$
such that
$$
[i_{-J}(X)] \cap [i_J(X)] = [0] \in H_2(Y, \Z).
$$
It can be illustreted as follows. We have an involution $\si$
$$
\si: Y \to Y
$$
which is represented by quaternionic real structures on the fibers
$$
[J] \in \pi^{-1}(pt) \mapsto [-J].
$$
Each fiber hasn't real points, and obviously
the intersection
$$
i_J(X) \cap i_{-J}(X) = \emptyset
$$
is trivial for each $J$.

One can compare this picture with the construction of $\Cal S$ above. Of course, the ambient spaces are the same, and 
one can reduce the problem from conformal classes to metrics by fixing
some appropriate nondegenerated positive (with respect to the orientation) 4- form and then for every conformal class
there is unique riemannian metric from this class which has this fixed 
4- form as volume form.
So one can represent a fiber of the projection $\pi_1$ as follows.
Let one has the twistor space $Y, \Cal J_g$ with almost complex structure
defined by our conformal class. Then
$$
\pi_1^{-1}(*_g) = \{ X' \subset Y | X'\quad {\text satisfies:} \}
$$

1. $X'$ is a smooth representation of the class $([H])^* \in H_4(Y, \Z)$

2. $X'$ intersects all the fibers transversally.

Due to this twistor picture one can prove the following
\proclaim{Proposition 1} Let $(g_0, J, \om) \in \sM_X^{K_i}$
has nontirivial image $\tau(g_0, J, \om) \in K^+$ so
$$
\tau(g_0, J, \om) \neq [0].
$$
Then the same is true for any $(g_0, J_1, \om_1)$.
\endproclaim

To prove this Proposition one can consider for a fixed metric $g_0$
the corresponding twistor space
$$
\pi: Y \to X
$$
with the almost complex structure $\Cal J_{g_0}$ defined by the conformal class
of $g_0$. But moreover using the standard metric of $\C\proj^1 = \pi^{-1}(p)$,
given metric $g_0$ and a connection on the $PU(2)$- bundle corresponding to
the Levi - Civita connection of the metric $g_0$ one gets as well
{\it a riemannian metric} $G_{g_0}$ on the whole $Y$ (it is much more clear   in the setup of the construction of $\Cal S$ (2.8) above ---  
all the necessary definitions were given there). So our riemannian
metric $g_0$ defines {\it a hermitian triple} $(G_{g_0}, \Cal J_{g_0}, \Om_{g_0})$ on $Y$ and we have the following
inclusion
$$
h: \Cal H_X \to \sM_Y
$$
and it's easy to see that the image lies in the following component 
$$Im h \subset \sM_Y^{K_{g_0}} \subset \sM_Y.
$$

Furhter, let us consider the second element $J$ of the original triple.
It corresponds to an inclusion of $X$ to $Y$ as described above. 

It's obvious that
$$
\aligned
d \pi (\Cal J_{g_0} &|_{i_{J}(X)}) = J\\
d \pi (G_{g_0} &|_{i_{J}(X)}) = g_0\\
d \pi (\Om_{g_0} &|_{i_{J}(X)}) = \om,\\
\endaligned
\tag 3.3
$$
and the triviality or nontrivialty of $\tau(g_0, J, \om)$
 can be related to the question of the triviality or nontriviality of the corresponding
class 
$$\wave{\tau}(G_{g_0}, \Cal J_{g_0}, \Om_{g_0}) \in H^2(Y, \R), \tag 3.4$$
restricted on the submanifold $i_{J}(X) \subset Y$.
Here the map 
$$\wave{\tau}: \sM_Y \to H^2(Y, \R)
$$
has the same definition as the map $\tau$ above (but really there is no
$K^+$ in $H^2(Y, \R)$ so the analogy isn't absolutely direct).

The relation between the following two integrals
$$
\int_X (\om)_H^2
\tag 3.4
$$
and
$$
\int_{i_J(X)} (\Om_{g_0})^2_H,
\tag 3.5
$$
where "H" denotes the corresponding harmonic components, defined by (0.1),
can be derived as follows. These two integral don't concide a priori, because of
the "vertical" component of the wedge square in the second one. But using our involution $\si$
one can "kill" this vertical component, namely
$$
\int_X (\om)_H^2 = \half \{ \int_{i_J(X)} (\Om_{g_0})^2_H -
\int_{\sigma(i_J(X))}(\Om_{g_0})^2_H \}.
\tag 3.6
$$
But since $\sigma(i_J(X)) = i_{-J}(X)$ and because of the harmonicity of $(\Om_{g_0})_H$ (3.6) can be reduced to the equality
$$
\int_X (\om)_H^2 = <[(\Om_{g_0})_H]; \pi^*(K_J)> \in \R,
\tag 3.6'
$$
where the right side {\it is a topological constant}.
Thus one can repeat all the discussion above to an arbitrary
$J_1$ with canonical class $K_{J_1}$, getting the same
equality (3.6'). This means that the expression
in the left side of (3.6') depends only on canonical class of our almost complex structure.

So the proof of the Proposition 1 is completed.

But one could note that we got even more, then the statement of 
Proposition 1. Namely the formula (3.6') ensures that
the riemannian norm of the corresponding harmonic component
$(\om)_H$ is constant under the contineous 
variations of the second element in our given hermitian triple.

Really, since all $(\om)_H$ are self dual, then
$$
\Vert (\om)_H \Vert_{g_0} =
\sqrt{\int_X (\om)_H^2} = const.
$$
Thus for a given metric $g_0$ we have a correspondence
$K_{J_i} \to \R_{\geq 0}$ between canonical classes and a set or non negative
real numbers.

\head \S 4. Semitwistor spaces
\endhead

Now we know that nontriviality of $\tau$ is stable with
respect to deformations of the second elements of the triples.
The next step is to establish the same stability
with respect to deformations of the first elements.

For this we'll study a new object: an anolog of the twistor space,
which depends on a fixed {\it almost complex structure} (instead of a riemannian metric).  Above we got
twistor space as total space for all (pointwise)  almost complex structures
compatible with a given riemannian metric. If one fixed vice versa an almost complex structure $J$ one can consider  all riemannian metrics pointwise compatible with $J$. Let us input some linear algebra to establish what is
a fiber in this case.

For $\R^4$ with given standard complex structure
$$
J = \left(\matrix 0& -E\\ E&0\\ \endmatrix \right)
$$
we are interested in such positive defined  bilinear forms which are preserved by $J$. Every form can be represented by a symmetric matrix
$A_Q$. Direct calculations give us the following answer
$$
A_Q = \left(\matrix a&b&0&0\\b&c&0&0&\\0&0&a&b\\0&0&b&c\\ 
\endmatrix \right),
$$
where $a, b, c \in \R$ and
$$
ac - b^2 > 0, \quad \quad a+c > 0,
\tag 4.1
$$
for positivity.

The set $\Cal L^+$ of points satisfied the two conditions above
is the inner half part of    quadratic cone $ac-b^2 =0$
 which belongs to the condition $a+c >0$.
Points on the boundary cone $ac - b^2 = 0$ correspond to {\it degenerated} metrics.

In turn to the conformal classes (of compatible riemannian metrics)
we get a bundle
$$
w: W \to X
\tag 4.2
$$
where $w^{-1}(pt)$ is topologically {\it a  two - dimensional disc}.
Let us call the total space of this bundle $W$ {\it semitwistor space}
corresponding to the given $J$.
A reason
to call it so is as follows: let us take  together with $\Cal L^+$
the other part
$$
\Cal L^- = \{(a, b, c)| \quad ac-b^2>0 \quad {\text but} \quad a+c <0 \}.
$$
This part corresponds to {\it negative defined} metrics, but nevertheless
after "projectivization" we get a sphere $S^2_p$ with marked equator
$S^1_p \subset S^2_p$ containing conformal classes of degenerated metrics.
Globalizing this picture over whole $X$ one gets a bundle on 2- spheres
$$
t: T \to X,
\tag 4.3
$$
where $t^{-1}(p) = S^2_p$, and this bundle, at least topologically,
corresponds to twistor bundle $\pi: \proj^1 \to X$  above.

But now we'd like to consider the following construction. Let us
realise the disc, which belongs to the set of all conformal classes, compatible
with a given almost complex structure $J$ as {\it the Lobachevsky plane}.
Really, on $\Cal L^+$ we have the natural norm
$$
v=(a, b, c) \mapsto \vert v \vert = ac-b^2
\tag 4.4
$$
and one can consider all vectors with the unit norm. They form
hyperbolic surface inside  $\Cal L^+$ which we'll denote
as $\Cal Lob$. It's clear that $\Cal Lob$ is naturally isomorphic to the Lobachevsky plane.
So one can regard
$$
w: W \to X
$$
as fibration on the Lobachevsky planes. So as on $\Cal Lob$ there exists
special riemannian metric one can define the universal conformal class $*_G$ on the whole $W$ in the same way as it was defined for the universal almost complex structure on twistor space. Unfortunately an anolog for Levi- Civita connection
doesn't exist. To avoid this one can choose a compatible riemannian metric
$g_0$ to identify all the fibers
For each $s \in W$ the corresponding to $g_0$ Levi- Civita connection defines
the decomposition
$$
TW_s = TX_{w(s)} \times T(\Cal Lob)_s,
$$
and on the second component we have standard riemannian metric  and on the first component one has
tautological conformal class on $X$ corresponding to  this point $s$. So we have
the universal conformal class $*_J$ on our semitwistor space $W$ and hence
one can define an almost complex structure on $W$. Namely, this universal conformal class defines over each point the orthogonal decomposition of
the tangent space so one can define the direct product of
two complex structure - the given on $TX$ and the standard on $T(\Cal Lob)$.

So for every almost complex structure $J$ on $X$, equipped with an additional tool - a compatible riemannian metric $g_0$, we defined  a 6- dimensional open manifold $W$ with boundary $S$ with
universal pair $*_J$ and $\Cal J_J$. But another problem is that
$*_J$ doesn't admit a nondegenerated extension to the boundary $S$ and the same
happens with our universal almost complex structure $\Cal J_J$.

Hence we have an analogy with the twistor case but now the situation is
much more complicated. First of all while a twistor space is compact
(and is smooth 6- dimensional manifold), a semitwistor space
is not compact and is represented by an open almost complex 6- dimensional
riemannian manifold. The boundary can be described
as $S^1$- bundle over $X$ and is denoted as 
$$
\partial W = S \to X.
$$

Now the description of all {\it globally defined} conformal classes
compatible with our given $J$ is the following. For each $*_g$ on based manifold $X$ one has
the corresponding {\it smooth inclusion}
$$
i_{*_g}(X) \to W
$$
such that $i_{*_g}(X)$ is a 4- submanifold which intersects the fibers transversally.
Again we have
$$
\aligned
d w(*_J &|_{i_{*_g}(X)}) = *_g\\
dw(\Cal J_J &|_{i_{*_g}(X)}) = J.\\
\endaligned
\tag 4.5
$$

Now we want to repeat the argument, which took place in the previous section
to establish nontriviality of our canonical map $\tau$ for all metrics
compatible with our given almost complex structure $J$. But in the present case
one has to use two additional arguments and one additional tool.

First of all the construction of semitwistor spaces above deals with {\it conformal classes} instead of riemannian metrics. So we have to fix
a nondegenerated 4- form on $X$, denoted as $d\mu$  (the volume form) and then get 
a one- to- one correspondence between conformal classes and riemannian metrics with the same volume form (so this fixing defines a section of fibration
$$
\wave{W} \to W, \quad i_{d \mu} (W) \hookrightarrow \wave{W}
\tag 4.6
$$
where $\wave{W}$ is the space of all riemannian metrics pointwise compatible with $J$).

Then one gets on $W$ {\it the universal riemannian metric} $G$ instead of the universal
conformal class $*_G$ and universal 2- form $\Om$. 

Further, the second point is that $W$ isn't compact and one has to impose
the following argument to establish stability of our canonical map $\tau$ with respect to changing of the first element in hermitian triple. Namely
instead of ordinary cohomology  one has to work with cohomology
with compact support. Then since our universal riemannian metric $G$
descends to zero on the boundary we can repeat all arguments on the Hodge decomposition of $\Om$ in terms of the cohomology with compact support and get 
the following statement:
\proclaim{Proposition 2} Let one fix a hermitian triple $(g, J, \om)$
with nontrivial image in $K^+ \subset H^2(X, \R)$. Then for every metric
$g_1$ compatible with $J$ and with the same volume form $d\mu$ as for the original one the image of $(g_1, J, \om_1)$ is nontrivial.
\endproclaim

An additional tool, which one needs, is an anolog of the involution
$\si:Y\to Y$, used in the previous section. Really, we again
have an involution of our semitwistor space $W$ due to the choosen
riemannian metric $g_0$. Namely, over each point of $X$ we have
the corresponding to $g_0$ framed point. So the involution can be defined as
the rotation of each fiber, centred in this framed points
$$
\si_{g_0}:W \to W; \quad \quad \si_{g_0}^2 = id.
$$

Then we can repeat with slight modifications all our arguments from the previous section. First of all, we have
$$
\int_X (\om_1)_H^2 = \half \{\int_{i_{g_1}(X)} (\Om_J)_H^2 +
\int_{\si_{g_0}(i_{g_1}(X))} (\Om_J)_H^2 \},
4.7
$$
where $\om_1$ corresponds to the fixed almost complex structure $J$
an a riemannian metric $g_1$, and the sign "+" in the right side of (4.7)
reflects the difference between "classical" involution
$\si$ for twistor spaces and "temporar" involution $\si_{g_0}$ (which has
the fixed point, for example, when $\si$ hasn't real points at all).

It remains just to use all tricks from the previous section. Namely,
one can rearrange (4.7) in "cohomological" style: for this
it's sufficient to note that since our fibers are simply connected
then $[i_{g_1}(X)]$ and $\sigma_{g_0} (i_{g_1}(X))$ are homologically the same
and since we've fixed the volume form these two submanifolds are sufficiently
"far" from the boundary.
Our harmonic 2- form $(\Om_J)_H$ descends to the boundary since our universal metric $G_J$ does, so it defines a compactly supported cohomological class
$[(\Om)_H] \in H^2_c(W, \R)$, thus
$$
\int_X(\om_1)_H^2 =  <[(\Om)_H]^2; [i_{g_1}(X)]> \in \R.
\tag 4.8
$$
These arguments can be exploited for any riemannian metric
with the fixed volume form, which gives us the statement of Proposition
2.

But again we get even more than this statement. Really,
we've got (with (4.8)) that for any two riemannian metric $g_0, g_1$, compatible with a given almost complex structure,  with the same
volume forms $d\mu_{g_0} = d \mu_{g_1}$ we have
$$
\Vert (\om_0)_H \Vert_{g_0} = \Vert (\om_1)_H \Vert_{g_1}.
$$
It is a marvel, because in this case we vary the norm as well as riemannian metric, but the value doesn't depend on such variation.

\head \S 5. From metric to conformal class
\endhead

As we've seen in the previous section the nontriviality condition
for $\tau$ is stable with respect to changing of the second element in 
hermitian triples. Now we'd like to consider the next step.
\proclaim{Proposition 3} Let for $(g, J, \om)$ the image
$\tau(g, J, \om) \in K^+_X$ is nontrivial. Then the same is true for any triple
of the shape $(e^f \cdot g, J, e^f \cdot \om)$ where
$f:X \to \R$ is a smooth real function 
$$
\tau(e^f \cdot g, J, e^f \cdot \om) \ne [0].
$$
 \endproclaim

\subheading{Proof} Let $(g, J, \om)$ is a triple with nontrivial
image $$\tau(g, J, \om) = [\om_H] \in H^2(X, \R) \quad |
\quad [\om_H] \ne [0]$$
where $\om_H$ is the corresponding harmonic 2- form with respect to
$*_g$ and the orientation. Consider the following integral
$$
\int_X e^f \om \wedge \om_H.
\tag 5.1
$$
Since both forms appearing in this integral are self dual this expression equals
to projection of the first form $e^f \om$ to the harmonic ray $\R <\om_H>$  modulo the norm of $\om_H$. So to prove the Proposition 3 it's sufficient to show that this integral isn't equal to zero.

For this   one can observe that if
$$
(\om, \om_H) \in C^{\infty}(X \to \R)
\tag 5.2
$$
is non negative real function then the integral above has to be strictly
positive.
So now our claim is to establish that the inner product $(\om, \om_H)$, induced
by our original riemannian metric $g$, is non negative everywhere. (And let's recall that $\om \wedge \om_H = (\om \wedge \om_H)$ everywhere because
of the self duality.)

It's clear that since
$$
\int_X \om \wedge \om_H = \int_X \om_H \wedge \om_H > 0
\tag 5.3 $$ (by the defintion), the function
$$
s = (\om, \om_H)
$$
has positive integral over $X$ with respect to the volume form $d\mu_g$.
So if there  exist points with negative volumes of $s$ then
it should impose the following picture. There is smooth 3- dimensional
submanifold $B \subset X$ which is 
$$
B = s^{-1}(0)
$$
and such that $X\setminus B = U^+ \cup U^-$ where 
$$
U^{\pm} = s^{-1}(\R^{\pm}).
$$

First of all, it's easy to see that $\om_H$ is equal to zero form on $B$
$$
\om_H |_B = 0.
$$
Really, let us consider smooth inclusion
$$
i: B \hookrightarrow X
$$
which gives us riemannian metric $i^* g$, nondegenerated everywhere
2- form $i^* \om$ and inherited orientation on $B$. So one can
use again the Hodge theorem over $B$ and it's clear that the harmonic
part of $i^*\om$ should be equal to $i^*\om_H$:
$$
i^* \om = i^* \om_H + d \rho_1 + d^* \rho_2
\tag 5.4
$$
{\it over $B$} where $\rho_1 \in \Om^1_B, \rho_2 \in \Om^3_B$.
Since we supposed that $(\om, \om_H) = 0$ over $B$, 
 we immediately get 
$$
\int_B (i^*\om_H, i^*\om_H)d\mu_{i^*g} = \int_B(i^*\om, i^*\om_H)d\mu_{i^*g} = 0
\tag 5.5
$$
so $B$ {\it has to be} a subset of zeroset for $\om_H$
$$
\om_H |_B = 0.
\tag 5.6
$$
These arguments can be used to establish what we need. Namely,
let $B_1$ is a smooth 3- dimensional submanifold of $X$ such that
$$B_1 \subset U^-$$
and
$$
i_1: B_1 \hookrightarrow X
$$
is the corresponding smooth inclusion. Then with respect to $i^*_1g$ and the inherited orientation over $B_1$ one has the same decomposition as (5.4)
$$
i^*_1 \om = i^*_1 \om_H + d \rho_3 + d^* \rho_4.
\tag 5.4'
$$
And the point is that
$$
\int_{B_1} (i^*_1 \om_H, i^*_1 \om_H) d \mu_{i^*_1 g} =
\int_{B_1}(i^*_1 \om, i^*_1 \om_H) d \mu_{i^*_1 g}
\tag 5.5'
$$
has to be negative so we get the contradiction.

From all above one can conlcude that:

1. $(\om, \om_H) \geq 0$ everywhere

and moreover

2. zerosets of $(\om, \om_H)$ and $(\om_H, \om_H)$ coincide.

Because of this for any smooth function $f \in C^{\infty}(X \to \R)$
we get
$$
\aligned
\int_X e^f \om \wedge \om_H = \int_X e^f (\om, \om_H) d \mu_g \geq\\
\int_X e^t (\om, \om_H) d\mu_g = e^t [\om_H]^2 > 0  \\
\endaligned
\tag 5.7
$$
where $t$ is the minimum of our smooth function $f$ on our compact smooth
manifold $X$. 

Hence we get the statement  of the Proposition 3.

\head \S 6. Good path, joining two triples
\endhead

Only one problem remains --- namely,
how one can join two triples from a component
$\sM^{K_i}_X \subset  \sM_X$ by an appropriate path, covered by  our Propositions 1, 2, 3  
which would give  us the statement
of the Main Theorem. Good path, joining
two triples $(g_0, J_0, \om_0)$ and
$(g_1, J_1, \om_1)$ can be constructed using the following
\proclaim{"Junction" Lemma}
Let $(g_0, J_0, \om_0) $ and $(g_1, J_1, \om_1)$ be a pair of hermitian triples
over $X$ with the same canonical class $K_{J_i} =
K \in H^2(X, \Z)$. Then there  exist
  an almost complex structure $J_{junct}$ with the same canonical class
together with two riemannian metrics $g_p, g_q$ such that

1) $g_p$ is
 compatible similtaneously with $J_0$ and $J_{junct}$

and

2) $g_q$ is compatible similtaniously with
$J_{junct}$ and $J_1$.
\endproclaim

So we'll construct a good path in $\sM_X^{K}$ joining two
given triple such that this path consists of five
"linear" pathes --- the first one is induced by changing
of riemannian metric with the fixed almost complex structure $J_0$,  the second is induced by changing of almost complex structures with the fixed riemannian metric $g_p$,  the third is induced by changing of riemannian metric with the fixed almost complex structure $J_{junct}$,  the fourth is induced by changing of almost complex structures  with the fixed riemannian metric $g_q$.
Of course, one has to impose the last "line" --- changing
of riemannian metrics from $g_q$ to $g_1$ with fixed
almost complex structure $J_1$.

The main step in this procedure is to construct
an appropriate almost complex structure $J_{junct}$ so
we begin 
 with the explanation of this point.

We work locally  over a point $x \in X$. So we have
in $\C\proj^3 = \proj( T^{\C}_x X)$
two pairs of projective lines $l_i, \Theta_x (l_i), i=0, 1,$
which correspond to our given almost complex structures $J_0, J_1$. It's clear that in general there is no a riemannian metric
(i.e.: a nondegenerated real quadric) in general which
contains all four projective lines. Let's consider what we have in $\C \proj^5 = \proj (\La^2 T^{*\C}_x X)$. There are
two pairs of points on the Grassmanian $Gr \subset \C \proj^5$ and our four projective lines lie on the same quadric in $\C \proj^3$ if and only if
the corresponding four points in $\C\proj^5$ lie on the same
real 2- plane (which is $\proj(\La^+_x)$ of course).

So given  these four points
$l_0, \Theta_x(l_0), l_1, \Theta_x (l_1)$ our aim is to construct such pair
of points $l_{junct}, \Theta_x(l_{junct})$ that:

the span $<l_i, \Theta_x(l_i), l_{junct}, \Theta_x(l_{junct})>$ is equal to projective 2- plane
$\proj^2_i, i=0, 1,$. The point is that
these two planes are automatically real (so
$\Theta_x(\proj^2_i) = \proj^2_i$) and one can
choose $g_p$ and $g_q$ such that $\proj^2_0 \subset
\C\proj^5$ is the projectivization of $\La^+_{g_p}$
and $\proj^2_1$ is the projectivization of $\La^+_{g_q}$.

So the remaining part of the proof is just an exercise 
in classical projective geometry. To construct
two conjugated points $l_{junct}, \Theta_x (l_{junct})$
one can use the following procedure. Let us take
the set of all projective 2-planes which contain the projective line $<l_0, l_1>$. Then if we find
such a 2- plane containing this line that it
has {\it two real points } $r_0, r_1$ (such points that
$$
\Theta_x (r_i) = r_i)
$$
our problem will be solved. Really, let's suppose that we've found such projective 2- plane $\pi$ that
$$
l_0, l_1, r_0, r_1 \in \pi
$$
where $l_0, l_1$ our given points and $r_0, r_1$
some distinct real points. Then  there exists 
conjugated 2- plane $\Theta_x (\pi)$
which contains the points $\Theta_x(l_0), \Theta_x(l_1),
r_0, r_1$. So one can construct the following two points
$$
\aligned
&l_{junct}' = <l_0, r_0> \cap <l_1, r_1> \in \pi \\
\Theta_x (&l_{junct}') =
<\Theta_x(l_0), r_0> \cap <\Theta_x(l_1), r_1> \in
\Theta_x (\pi).\\
\endaligned
\tag 6.1
$$
It's easy to see that projective lines $<l_0, \Theta_x(l_0)>$
and $<l_{junct}', \Theta_x(l_{junct}')>$ lie in the same
projective 2- plane (which {\it isn't our} $\pi$!)
and the same is true for projective lines $<l_1, \Theta_x(l_1)>$ and $<l_{junct}', \Theta_x(l_{junct}')>$.

But two constructed points aren't what we need because
a priori these points don't lie on our Grassmanian $Gr$
so don't represent any lines in our $\C\proj^3$. 
So one has to construct the intersection
$$
<l_{junct}', \Theta_x (l_{junct}')> \cap Gr
= \{ l_{junct}, \Theta_x (l_{junct}) \}
\tag 6.2
$$
and these two points would be our
$l_{junct}$ and $\Theta_x(l_{junct})$. It's not 
hard to see that this intersection is a pair of points
indeed because our quadric $Gr \subset \C\proj^5$ is real
with respect to $\Theta_x$ and our projective line
$<l_{junct}', \Theta_x(l_{junct}')>$ is real but without
real points (by the construction).

So it remains to explain why there exists such appropriate
projective 2- plane $\pi$ which contains given two points
$l_0, l_1$ and a pair of distinct real points. 

We construct $\pi$ as follows:

1. consider the span $$<l_0, l_1, \Theta_x(l_0),
\Theta_x(l_1)> = \C\proj^3 \subset \C\proj^5; \tag 6.3$$

2. our standard real structure $\Theta_x$ restricted
on $\C\proj^3$ gives us {\it a standard real structure} (which we denote as $\Theta_x'$)
on $\C \proj^3$;

3. since $\Theta_x'$ is a standard real structure
there exists invariant with respect to $\Theta_x'$ subset
$\R\proj^3 \subset \C\proj^3$;

4. consider a 2- plane $\pi$, containing the projective line
$<l_0, l_1>$ and lieing in our $\C\proj^3$;

5. since every two projective planes in $\C\proj^3$ 
have non trivial intersection (which is a projective line)
we get the projective line
$$
\rho = \pi \cap \Theta_x'(\pi);
\tag 6.4
$$

6. this projective line $\rho$ is invariant with respect to
the $\Theta_x'$- action (i.e. is a real projective line)
$$
\Theta_x'(\rho) = \rho;
$$

7. every real projective line in $\C\proj^3$
has at least two real points.

One can see that really {\it every } 2- plane,
contained in the span (6.3), satisfies the required property (i.e. has two distinct real points).

We have to check only step 2 and step 7 from the list above.
The other steps are rather obvious.

First of all, let us recall that there exist two types of
real structure on projective spaces. The first one is standard real structure and in appropriate homogeneous 
coordinates in $\C\proj^n$ it has the  form
$$
\Theta_s (z_0: ...: z_n) =
(\overline{z}_0: ...: \overline{z}_n).
$$
 Evedentely, it has real points. The second type exists
only on $\C\proj^n$ with $n$ odd. In this odd dimensional case
we have quaternionic real structure:
$$
\Theta_q (z_0: z_1: ...: z_{n-1}: z_n) =
(-\overline{z}_1: \overline{z}_0: ...: -\overline{z}_n:
\overline{z}_{n-1}).
$$
This quaternionic structure has no real point (it can be checked directly). 

\subheading{Step 2} Since our $\C\proj^3$ is preserved by
$\Theta_x$ (the four points are transformed ones to the others) we get well defined real structure $\Theta_x'$
on $\C\proj^3$. We have to establish that this restricted
$\Theta_x'$ is a standard
real structure. For this it is sufficient to show
that there exists a real point (because there are only two possibilities --- standard and quaternionic). So consider
the projective line $<l_0, \Theta_x(l_0)> \subset \C\proj^3
\subset \C\proj^5$. This line is preserved by the real structure $\Theta_x$. Hence this reduction of the step 2 together with the step 7 will be proved if  we prove the following 
\proclaim{Statement} Let a projective line $\rho \in \C\proj^n$
is preserved by a fixed standard real structure
$\Theta_s$. Then there exist at least two real points
on this projective line.
\endproclaim

Proof. Let us fix appropriate homogeneous coordinates
on $\C\proj^n$ such that our given real structure
$\Theta_s$ has the form (??) in this coordinate system.
Then consider a system of homogeneous equations in the same coordinates, corresponding to our projective line $\rho$:
$$
\sum_{i=0}^n \al^i_j z_i = 0,
\tag 6.5
$$
$\al_j^i \in \C, j = 1, ..., n-1$. Because of the $\Theta_s$- invariance we get
that the points of $\rho$ satisfy the following system
$$
\sum_{i=0}^{n} \al^i_j \overline{z}_i = 0
\tag 6.5'
$$
as well as the original one. Then it is easy to see
that the points 
$$
r_0 = (Re z_0: ...: Re z_0), \quad r_1 =
(Im z_0: ... : Im z_0)
$$
should satisfy both systems. So we've found two real points.

Of course, there exist not only two real points on $\rho$.
There is a real projective line
$$
\rho_{\R} \cong \R\proj^1, \quad \rho_{\R} \subset \rho
$$
containing real points, but for our aim it is sufficient to
find only two distinct real points.

With the Statement in hand we prove the step 2 and the step
7, applying this result to the case $n=5$ and the case $n=3$
respectively.

It remains to note that one can globalize this local picture
over whole $X$ (using putching arguments) and then deduce
that such good path always exists.

So now we have a good "partially linear" path
joinig two given hermitian triples, which were choosen
arbitrary, and then applying Propositions
1, 2 and 3 along this path one gets the statement of
our Main Theorem.

\head \S 7. Additional remarks
\endhead

First of all it's necessury to say that the Main Theorem above gives us
new invariants of smooth structures in 4- geometry.

Really, it's clear that for two diffeomorphic 4- manifolds $X_1, X_2$
with some diffeomorphism
$$
\phi: X_1 \to X_2
$$
if $K_2 \in H^2(X_2, \Z)$ has {\it nontrivial generalized image}
in $H^2(X_2, \R)$ (so if for every triple $(g, J, \om)$ with
$K_J = K_2$ the image $\tau (g, J, \om) \ne [0] $) the  corresponding
class $\phi^* K_2 \in H^2(X_1, \Z)$ should be "nontrivial" in this sense
too.

Let us define the following subset in $H^2(X, \Z)$ for a smooth compact orientable 4- manifold $X$. A class $K_i \in H^2(X, \Z)$ is called
{\it PS}- class if the generalized image is nontrivial.
The subset $N_{PS} \subset H^2(X, \Z)$  {\it is preserved by the group} $Diff^+ X$ (as the basic classes in the gauge theories).
Of course, the generalization of the Taubes result ([9]) gives us the following inclusion
$$
N_{PS} \subset \{ {\text the}\quad {\text set}\quad
{\text of}\quad {\text basic}\quad {\text classes} \}
$$
but we'd like to emphercise that our definition of the PS- classes
works in the case $b^+_2 = 1$ as well as in the "classical" for the gauge theories case $b^+_2 > 1$ (of course, if $b^+_2 = 0$ the subset
$N_{PS} $ is empty).

Moreover, with remarks from the ends of Sections 3 and 4 in hand one can prove the following
\proclaim{Theorem} Let $(g_1, J_1, \om_1)$ and $(g_2, J_2, \om_2)$
a pair of hermitian triples over $X$ s.t.
$$
K_{J_1} = K_{J_2}, \quad \quad d\mu_{g_1} = d \mu_{g_2}.
$$
Then 
$$
\Vert (\om_1)_H \Vert_{g_1} = \Vert (\om_2)_H \Vert_{g_2}.
$$
\endproclaim

This allows us to define a version of the PS- invariants, defined above.
Namely, consider the following real nonegative function
$$
\tau_{\R}: \Cal M \to \R_{\geq 0}
$$
s.t.
$$
(g, J, \om) \mapsto \Vert (\om)_H \Vert_g.
$$
Obviously this function is invariant under the $Diff^+(X)$- action.

But one could try to derive some more think invariants
of smooth structures from the construction of the canonical map $\tau$. Let us consider the complete image
$Im \tau \subset K^+ \subset H^2(X, \R)$. It's clear again
that 

1. $Im \tau$ is a "subcone" in $K^+$. The point is
that if some point $k \in K^+$ is realized as the image of 
a triple $(g, J, \om)$ then for every positive real number
$r>0$ the point $r \cdot K \in K^+$ is realized too;

2. "subcone" $Im \tau$ is a connected subset and the canonical map $\tau$ is a contineous map;  

3. it's quite natural to projectivize this "subcone"
and then the topology of the porjective "manifold"
reflects properties of the underling smooth structure.
Namely if there are two smooth structures
on the same topological based manifold
then two corresponding "projectivized subcones"
should be topologically ismorphic one to other.

But now it's the time to finish this article hoping
that the new notion of pseudo symplectic manifolds
inspired by the considerations in the framework of Seiberg --- Witten theory will lead to new results in the smooth
classification of 4- manifolds.

\head Asknoledgements
\endhead

 I'd like to express my gratefull to
all stuff of Max- Planck- Institute
in Bonn but especially to Prof. Yu. I. Manin
and Prof. G. Faltings for the great possiblity
to work in the atmosphere of friendel attention
for usefull discussions and remarks and for hospitality at all.

Also I'd like to thank Prof. V. Shokurov (JHU, Baltimor), Prof.
F.Bogomolov and Prof. S.Capell (CU, NYU) for their attension,
valuable discussions and remarks, which inspired the cleaning
of the previous version of this text.

\enddocument